\documentclass{amsart}
\usepackage[centering]{geometry}

\usepackage{amsfonts}
\usepackage{amssymb}
\usepackage{amsmath}
\usepackage{hyperref}
\usepackage{graphicx}
\usepackage{color}

\newcommand{\half}{\frac{1}{2}}
\newcommand{\thalf}{\tfrac{1}{2}}
\newcommand{\summ}{\mathop{{\sum}^{\star}}}

\numberwithin{equation}{section}

\newtheorem{theorem}{Theorem}[section]
\newtheorem*{theorem*}{Remark}
\newtheorem{proposition}[theorem]{Proposition}
\newtheorem{lemma}[theorem]{Lemma}

\makeatletter
\@namedef{subjclassname@2010}{%
  \textup{2010} Mathematics Subject Classification}
\makeatother

\begin{document}

\title{Non-vanishing of Dirichlet $L$-functions in Galois orbits}

\author{Rizwanur Khan}
\address{
Science Program\\ Texas A\&M University at Qatar\\ PO Box 23874\\ Doha, Qatar}
\email{rizwanur.khan@qatar.tamu.edu, hieu.ngo@qatar.tamu.edu }

\author{Djordje Mili\'cevi\'c}
\address{ Department of Mathematics\\ Bryn Mawr College\\ 101 North Merion Avenue\\ Bryn Mawr, PA 19010\\ U.S.A}
\email{dmilicevic@brynmawr.edu}

\author{Hieu T. Ngo}

\thanks{D.M. acknowledges support by the National Security Agency. Project is sponsored by the NSA under Grant Number H98230-14-1-0139. The United States Government is authorized to reproduce and distribute reprints notwithstanding any copyright notation herein.}

\begin{abstract} 
A well known result of Iwaniec and Sarnak states that for at least one third of the primitive Dirichlet characters to a large modulus $q$, the associated $L$-functions do not vanish at the central point. When $q$ is a large power of a fixed prime, we prove the same proportion already among the primitive characters of any given order. The set of primitive characters modulo $q$ of a given order can be described as an orbit under the action of the Galois group of the corresponding cyclotomic field. We also prove a positive proportion of nonvanishing within substantially shorter orbits generated by intermediate Galois groups as soon as they are larger than roughly the square-root of the prime-power conductor.
\end{abstract}

\maketitle

\section{Introduction}

Central values of $L$-functions are of fundamental importance in number theory. In particular, a host of results and conjectures, including the Birch and Swinnerton-Dyer Conjecture, the Riemann Hypothesis, and the Katz-Sarnak Density Conjecture, predict in various contexts that the central values of $L$-functions (or their derivatives, as appropriate for root number reasons) hold key arithmetic information and should vanish only when there are deep arithmetic reasons for them to do so and that this should be an exceptional occurrence in suitably generic families.

Introduced by Bohr and Landau \cite{bohr-landau} in their study of zeroes of the Riemann zeta-function and notably used by Selberg \cite{selberg} in the course of proving that a positive proportion of these zeroes lie on the critical line, the ``mollifier'' is the most versatile tool used in analytic number theory to prove the non-vanishing of central values of $L$-functions in families, often achieving a positive proportion result. What is by now a classical result using the mollifier is one of Iwaniec and Sarnak \cite{iwasar}, concerning non-vanishing in the family of Dirichlet $L$-functions.  They proved that for at least $(\frac{1}{3} -\epsilon)$ of the primitive Dirichlet characters modulo $q$, where $q$ is any integer sufficiently large in terms of $\epsilon$, the central value $L(\frac{1}{2},\chi)$ is not zero. This is currently the best known result that can be proved with a ``one-piece'' mollifer. Earlier, Balasubramanian and Murty \cite{balmur} had established a smaller positive proportion of non-vanishing, and recently, Bui \cite{bui} proved, using a``two-piece'' mollifier, that about $34\%$ of the central values in this family do not vanish. When one restricts to the quadratic Dirichlet $L$-functions, Soundararajan~\cite{sou} established that for at least $\frac{7}{8}$ of the fundamental discriminants $|d|\le X$, the central value $L(\frac{1}{2}, (\frac{d}{\cdot}))$ is not zero, as $X\to \infty$. It is generally conjectured (see the discussion in \cite{sou}) that \begin{equation}
\label{nonvanish}
L(\tfrac12, \chi)\neq 0
\end{equation}
for any primitive character $\chi$. Such a statement appears to be substantially beyond the reach of currently available technology.

Let $\xi$ be a primitive $\phi(q)$-th root of unity, where $\phi$ is the Euler totient function. The Galois group $G=\text{Gal}(\mathbb{Q}(\xi)/\mathbb{Q})$ acts on the set of primitive Dirichlet characters modulo $q$ as follows. For $\sigma\in G$, we define $\chi^{\sigma}$ to be the character given by $\chi^\sigma(n)=\sigma(\chi(n))$ for all integers $n$. The Galois action partitions the set of characters into orbits, which are particularly natural from the arithmetic point of view, since the associated Dirichlet $L$-functions (by definition) share the same field of coefficients. One is led to wonder whether a positive proportion of non-vanishing can be proven for Dirichlet $L$-functions within each orbit. The aforementioned results \cite{iwasar} and \cite{bui} do not preclude the possibility that $L(\half, \chi)=0$ for all $\chi$ in some orbit $\mathcal{O}$. 

As a positive proportion statement toward \eqref{nonvanish} in this context, we conjecture that for some $c>0$, we have that
\begin{align}
\label{conj} \frac{1}{|\mathcal{O}|}\sum_{\substack{\chi\in \mathcal{O}\\ L(\half, \chi)\neq 0}} 1 \ge c -\epsilon
\end{align}
for any $\epsilon >0$, any orbit $\mathcal{O}$ of cardinality $|\mathcal{O}|>q^{\epsilon}$, and any integer $q$ sufficiently large in terms of $\epsilon$. One might hope to match the best constant $c=0.34$ currently available for the full set of primitive Dirichlet $L$-functions or the more classical proportion of Iwaniec and Sarnak. We establish the latter in the case that $q$ is a large power of a prime. 
\begin{theorem} \label{mainthm} Let $q=p^k$ for an odd prime $p$. For any $\epsilon>0$ and $k$ large enough in terms of $\epsilon$ and $p$, we have that (\ref{conj}) holds with $c=\frac{1}{3}$.
\end{theorem}
When $q=p^k$, the orbits under the Galois action can be described as follows. First recall that a necessary and sufficient condition for the existence of a \emph{primitive} character modulo $p^k$ having order $l$ is that $l=p^{k-1}d$ for some $d\mid (p-1)$. The set of all primitive Dirichlet characters modulo $q$ whose orders equal $p^{k-1}d$ forms an orbit $\mathcal{O}$ of cardinality
\begin{align}
\label{size} |\mathcal{O}|=\phi(p^{k-1}d),
\end{align}
and every orbit arises in this way. Thus $\mathcal{O}$ depends on $d$ but we suppress this in the notation. For these facts, see, for example, \cite[Chapter 5]{dav} and \cite[page 16]{chi}.

Although, from an analytic perspective, it may appear that we are dealing with a family only slightly thinner than the original unitary family of all primitive Dirichlet characters, in reality this gives rise to a significant difficulty, which we describe below in the introduction. The key feature here is that the family is ``thinning out'' in a thoroughly arithmetic (rather than analytic) way. In fact, the same device that we use to overcome this basic difficulty subsequently allows us to prove a positive proportion of non-vanishing in substantially smaller ``thin orbits'' in Theorem~\ref{thinorbitthm} below.

\bigskip

The problem of studying the non-vanishing of $L$-functions within Galois orbits is a natural one that has yielded some of the strongest known results in the subject. Let $q=p^k$ for the remainder of this paper. Let $f$ be a holomorphic newform of weight 2 and level coprime to $p$ which has rational Fourier coefficients (equivalently, $f$ is associated to an elliptic curve over $\mathbb{Q}$ of conductor coprime to $p$, by \cite[Theorems 7.14, 7.15]{shi3} and \cite{wil, taywil, bcdt}). Let $L(s, f \times \chi)$ be the $L$-function of $f$ twisted by a primitive Dirichlet character $\chi$ modulo $q$, having central point $s=1$ by a functional equation normalized so as to relate $L(s,  f \times \chi)$ and $L(2-s, f \times \overline{\chi})$. Rohrlich \cite{roh} showed that $L(1,f\times \chi)$ does not vanish as long as $k$ is large enough in terms of $p$ and $f$. To prove his result, Rohrlich appealed to an ``algebraicity'' theorem of Shimura \cite[Theorem 1]{shi1,shi2}, which implies that if $L(1, f \times \chi)=0$ then $L(1, f \times \chi^\sigma)=0$ for all $\sigma \in G$. By this, if the sum
\begin{align}
\label{roh-sum}  \sum_{\chi \in \mathcal{O}} L(1, f \times \chi)
\end{align}
is nonzero, then every summand is nonzero. Rohrlich found an asymptotic for (\ref{roh-sum}) for any orbit $\mathcal{O}$ when $k$ is large enough in terms of $p$ and $f$, and showed that the main term is indeed nonzero. Chinta \cite{chi} extended Rohrlich's work to the case of prime moduli $q$ by considering instead the sum
\begin{align}
\label{chi-sum} \sum_{\chi \in \mathcal{O}} L(1, f \times \chi) M(f\times \chi),
\end{align}
where $M(f\times \chi)$ is a truncation of the formal Dirichlet series for $L(1, f \times \chi)^{-1}$. The effect of this mollifier is that each summand of (\ref{chi-sum}) is ``morally'' close to 1, and this allows Chinta to show that the sum is nonzero provided $|\mathcal{O}|>q^{\frac{7}{8}+\epsilon}$ and $q$ is large enough in terms of $f$. In particular, this implies the non-vanishing of $L(1,f\times \chi)$ over big orbits when $q$ is a large enough prime in terms of $f$, a fact which does not follow from Rohrlich's result.

All these results for twists of elliptic modular $L$-functions rely heavily on the algebraicity results of Shimura. Such a route is not available in the present context of central values of Dirichlet $L$-functions.

\bigskip

We establish Theorem \ref{mainthm} by evaluating the mollified moments
\begin{align}
\label{mol1} \sum_{\chi \in \mathcal{O}} L(\thalf, \chi) M(\chi)
\end{align}
and
\begin{align}
\label{mol2} \sum_{\chi \in \mathcal{O}} |L(\thalf, \chi)|^2 |M(\chi)|^2,
\end{align}
where
\begin{align}
\label{moldef} M(\chi)=\sum_{m\le q^\theta} \frac{a_m \chi(m)}{m^{\half}}
\end{align}
is a mollifier of length $q^\theta$, for some $\theta\ge 0$ and coefficients $a_m$ satisfying $a_m\ll m^{\epsilon}$ and $a_1=1$. We are able to obtain asymptotics with a power-saving error term for arbitrary mollifiers when $\theta$ is any fixed constant satisfying $0\le \theta < \frac{1}{2}$; see Propositions~\ref{1stmol} and \ref{2ndmol}. As discussed in section~\ref{MollifierSubsection}, this allows us to deduce Theorem~\ref{mainthm} with the proportion of non-vanishing $c=\frac13$ by taking $\theta\to\frac12$.

Although $|L(\half, \chi)|^2$ can be considered to be analogous to $L(1,f \times \chi)$, our problem has some important differences from the one considered by Rohrlich and Chinta. Firstly, the additional factor $|M(\chi)|^2$ in (\ref{mol2}) makes our problem more complex. This is readily seen when comparing with (\ref{roh-sum}), and to compare with (\ref{chi-sum}) we note that Chinta's method only works for a specific mollifier while ours works for an arbitrary mollifier, and that $M(f\times \chi)$ has length at most $q^{\frac{1}{4}-\epsilon}$ (see \cite[pg 22]{chi}) while $|M(\chi)|^2$ has length at most $q^{1-\epsilon}$. However one must keep in mind that Chinta proves a result which works for $q$ prime, while ours does not.

The second difference is the method of proof. To evaluate (\ref{roh-sum}), Rohrlich had to consider the averages
\begin{align}
\label{rohavg}  \sum_{\chi\in \mathcal{O}}  \chi(n)
\end{align}
for $n\le q^{1+\epsilon}$, while for (\ref{mol2}) we must consider the averages
\begin{align}
 \label{ouravg} \sum_{\chi\in \mathcal{O}}  \chi(n_1m_1)\overline{\chi}(n_2m_2)
\end{align}
in which $n_1n_2\le q^{1+\epsilon}$ and $m_1,m_2\le q^{\theta}$. As we will see, (\ref{rohavg}) is zero unless $n^{p-1} \equiv 1 \bmod p^{k-1}$. From this, Rohrlich could immediately conclude that $n=1$ or $n> p^{\frac{k-1}{p-1}}$, thereby effectively isolating the contribution of the diagonal term $n=1$. In contrast,  (\ref{ouravg}) is zero unless
\begin{align}
\label{ourcondition} (n_1m_1)^{p-1}\equiv (n_2m_2)^{p-1} \bmod p^{k-1}.
\end{align}
Writing $n_1m_1\equiv\zeta n_2m_2\bmod{p^k}$, it is now much harder to isolate the contribution of the diagonal terms $n_1m_1=n_2m_2$ (that is, $\zeta\equiv 1\bmod{p^k}$); \emph{a priori} it is, for example, perfectly plausible that $n_1m_1$ and $n_2m_2$ could be fairly close to each other without actually being equal. To isolate the contribution of the diagonal terms $n_1m_1=n_2m_2$, we will appeal to the $p$-adic version of Roth's theorem, Lemma~\ref{rohrlich-roth} below. The upshot is that, keeping in mind that $\zeta^{p-1}\equiv 1\bmod{p^k}$, having two solutions to \eqref{ourcondition} within the same class of $\zeta\not\equiv\pm 1\bmod{p^k}$ too close to each other would ultimately yield too good of an approximation in the $p$-adic norm to a $(p-1)^{\text{th}}$ $p$-adic root of unity.

It is interesting that Roth's theorem, a deep result from diophantine approximation, should be used to prove the non-vanishing of $L$-functions. This connection has been made before in different contexts by Rohrlich \cite{roh2} and Greenberg \cite{gre}. Our paper offers another such example and it seems to be the first one involving a family of Dirichlet characters as well as the first one involving a genuinely second mollified moment.

\bigskip

Taking $\theta=0$ in our evaluation of (\ref{mol1}) and (\ref{mol2}) yields in particular the first and second moments of $L(\half, \chi)$ over Galois orbits. Before stating this result, we note that the value $\chi(-1)$ is the same for every character in any given orbit $\mathcal{O}$, as it is a rational number.
\begin{theorem} \label{2ndthm} Let $q=p^k$ for an odd prime $p$ and let $\mathcal{O}$ be any Galois orbit of primitive Dirichlet characters mod $q$. Suppose that $\chi(-1)=(-1)^\iota$ for any $\chi\in \mathcal{O}$. We have that
\begin{align*}
&\frac{1}{|\mathcal{O}|} \sum_{\chi \in \mathcal{O}} L(\thalf, \chi)=  1 + O\big(q^{-\frac14+\epsilon}  \big) \\
&\frac{1}{|\mathcal{O}|}  \sum_{\chi \in \mathcal{O}} |L(\thalf, \chi)|^2=  \frac{p-1}{p} \Big(  \log\Big(\frac{q}{\pi }\Big) +\frac{\Gamma'(\frac{1+2\iota}{4})}{\Gamma(\frac{1+2\iota}{4})} +2\gamma+2\frac{\log p}{p-1} \Big) + O\big(q^{-\frac14+\epsilon}  \big) ,
\end{align*}
for any $\epsilon>0$, where the implied constants depend on $\epsilon$ and $p$, and $\gamma=0.57721\cdots$ is the Euler constant.
\end{theorem}

We remark that the implicit constants in all our main results are ineffective due to their dependence on $p$-adic Roth's Theorem. However, the first moment in Theorem 1.2 and the mollified first moment in Proposition~\ref{1stmol} can also be evaluated without recourse to $p$-adic Roth's Theorem, at the expense of the error terms $O(q^{-\frac14+\epsilon})$ being replaced by the weaker but effective error terms $O(q^{-\frac1{2(p-1)}+\epsilon})$. This will be shown in the course of the proof of Proposition~\ref{1stmol}.

\bigskip

Finally, we address the refined question of non-vanishing in smaller sub-families within the Galois orbits of primitive characters modulo $q$. A rather natural sub-family emerges when considering orbits of primitive characters under various subgroups $H$ of the Galois group $G=\text{Gal}(\mathbb{Q}(\xi)/\mathbb{Q})$. These subgroups form a partially ordered set (corresponding by Galois theory to the tree of intermediary field extensions $\mathbb{Q}\leq K\leq\mathbb{Q}(\xi)$), and, as the subgroup $H$ varies from $G$ through its various subgroups to the identity, the corresponding orbits of a fixed primitive character modulo $q$ can be seen as interpolating (or shrinking) between its full Galois orbit, considered in Theorem~\ref{mainthm}, and the individual character.

We describe these ``thin orbits'' explicitly in cases of our interest. For every $0\leq\kappa\leq k-1$, denote $K_{k-1-\kappa}=\mathbb{Q}\big(\xi^{p^{\kappa}}\big)$. (Note that the field $K_{\ell}$ is independent of $k$.) Since $[K_{k-1}:K_0]=\phi(p^{k-1})\asymp_p\phi(\phi(q))$ and we are primarily concerned with the case of fixed $p$ and large $k$, we focus here on the tower of these intermediate fields
\[ \mathbb{Q}(\xi)=K_{k-1}\supseteq\dots\supseteq K_0\supseteq\mathbb{Q}. \]
The Galois group ${\rm Gal}(\mathbb{Q}(\xi)/\mathbb{Q})$ acts transitively on any given Galois orbit $\mathcal{O}$. The intermediate Galois group ${\rm Gal}(\mathbb{Q}(\xi)/K_{k-1-\kappa})$ therefore acts on $\mathcal{O}$; we call an orbit of this action a \emph{thin Galois orbit}, and we write $\mathcal{O}_{\kappa}$ for any one of these thin orbits. Note that already the thin orbits $\mathcal{O}_{k-1}$ refine the full Galois orbits $\mathcal{O}$, with thin orbits $\mathcal{O}_{\kappa}$ for smaller $\kappa$ being progressively smaller (so that we may think of the parameter $\kappa$ essentially as an indicator or the logarithmic size of the corresponding thin orbits), all the way to the extreme case of $\kappa=0$, which corresponds to the single primitive characters.

It is not difficult to see that $\sigma\in\mathrm{Gal}(\mathbb{Q}(\xi)/\mathbb{Q})$ satisfies $\sigma\in\mathrm{Gal}(\mathbb{Q}(\xi)/K_{k-1-\kappa})$ if and only if
\begin{equation}
\label{SigmaXiCondition}
\sigma(\xi)=\xi^a\quad\text{for some }a\equiv 1\bmod{p^{k-1-\kappa}(p-1)}.
\end{equation}
We thus see that, equivalently, two characters $\chi_1,\chi_2$ are in the same thin Galois orbit $\mathcal{O}_{\kappa}$ if and only if $\chi_1\overline{\chi_2}$ has order dividing $p^{\kappa}$ (informally speaking, if $\chi_1$ and $\chi_2$ differ by an ``algebraically simpler'' character, one that is ``shallower'' in the sense of the ``depth aspect'' of modulus $q=p^k$ with large $k$). In particular, all thin orbits $\mathcal{O}_{\kappa}$ for the same $\kappa$ are of equal size given by (keeping in mind that $(a,\phi(q))=1$)
\[ |\mathcal{O}_{\kappa}|=\begin{cases} p^{\kappa}, &0\leqslant\kappa<k-1,\\ \phi(p^{k-1}), &\kappa=k-1.\end{cases} \]

\def\im{\mathop{\mathrm{im}}}

Another, more ``analytic'', way to think about the thin orbits $\mathcal{O}_{\kappa}$ is provided by the explicit characterization of the ``principal part'' of the dual of the group $(\mathbb{Z}/q\mathbb{Z})^{\times}$ for a high prime power $q=p^k$, which is essentially due to Postnikov~\cite{Postnikov}. Let $X_k=(\mathbb{Z}/p^k\mathbb{Z})^{\times}$, and let $\pi_a$ denote the power map $\pi_a:X_k\to X_k$, $[x]\mapsto [x^a]$. Corresponding to the decomposition $X_k=X_{k0}\times X_{k1}$, where $X_{k0}=\im\pi_{p^{k-1}}=\ker\pi_{p-1}$, $|X_{k0}|=p-1$, and $X_{k1}=\im\pi_{p-1}=\ker\pi_{p^{k-1}}=\{[x]\in X_k:x\equiv 1\bmod p\}$, we have the canonical decomposition of dual groups $\hat{X}_k\cong\hat{X}_{k0}\times\hat{X}_{k1}$. Let $\log_p$ denote the $p$-adic logarithm and $\psi(x)$ denote the ``standard'' additive character $\psi:\mathbb{Q}_p\to\mathbb{C}^{\times}$ such that its kernel is exactly $\mathbb{Z}_p$ and that $\psi(x)=e^{2\pi ix}$ for $x\in\mathbb{Z}[1/p]\subseteq\mathbb{Q}_p\cap\mathbb{R}$. According to Postnikov's lemma (for $p>2$; see also \cite[Lemma~13]{MilicevicSubWeyl}), every character $\chi^{(1)}\in\hat{X}_{k1}$ is of the form
\[ \chi^{(1)}(1+pt) = \chi^{(1)}_a(1+pt)=\psi\left(\frac{a_0\log_p(1+pt)}{p^k}\right) \]
for some $a=a_0p^{-k}\in p^{-k}\mathbb{Z}_p/p^{-1}\mathbb{Z}_p$, with primitive characters corresponding to $a\in p^{-k}\mathbb{Z}_p^{\times}/p^{-1}\mathbb{Z}_p$.

The isomorphism of $p^{-k}\mathbb{Z}_p/p^{-1}\mathbb{Z}_p\to\hat{X}_{k1}$ given by $a\mapsto\chi^{(1)}_a$ induces a metric on $\hat{X}_{k1}$ via $d(\chi^{(1)}_a,\chi^{(1)}_b)=|a-b|_p/p=\mathrm{cond}(\bar{\chi}^{(1)}_a\chi^{(1)}_b)$ for $\chi^{(1)}_{a}\neq\chi^{(1)}_b$. Relative to the above decomposition of $\hat{X}_k$ and this metric on $\hat{X}_{k1}$, the thin orbit $\mathcal{O}_{\kappa}$ containing a character $\chi=\chi^{(0)}\chi^{(1)}$ is, for $\kappa>0$, precisely the set $\{\chi^{(0)}\}\times B[\chi^{(1)},p^{\kappa}]$; here $B[\chi^{(1)},p^{\kappa}]$ denotes the closed ball in $\hat{X}_{k1}$ with center $\chi^{(1)}$ and radius $p^\kappa$ with respect to the above-defined metric. In particular, all characters in a thin orbit $\mathcal{O}_{\kappa}$ share the same $\hat{X}_{k0}$-component and their $\hat{X}_{k1}$-components are all close to each other, with the corresponding neighborhood around a fixed character $\chi$ shrinking as $\kappa$ decreases.

From the point of view of harmonic analysis, we see clearly the basic difficulty of isolating individual $\chi^{(0)}\in\hat{X}_{k0}$ in our orbits (which, in a modified form, is already present in isolating the full orbits $\mathcal{O}$), which on the dual side is reflected by the initial survival of the $(p-1)^{\text{th}}$ roots of unity in Lemma~\ref{chiavgthin}, followed by isolating characters $\chi^{(1)}$ in smaller neighborhoods within $\hat{X}_{k1}$, which corresponds to the survival of further terms in more permissive congruence classes containing these roots of unity.

\bigskip

Our techniques, which ultimately rely on the impossibility of overly good $p$-adic approximations to algebraic integers, are very well suited to the study of thin orbits of primitive characters to prime power moduli and give the following refinement of Theorem~\ref{mainthm}.

\begin{theorem} \label{thinorbitthm} Let $q=p^k$ for an odd prime $p$. For any $\epsilon>0$ and $k$ large enough in terms of $\epsilon$ and $p$, and for any $\kappa>k/2$, we have that (\ref{conj}) holds also when $\mathcal{O}$ is replaced by any ``thin orbit'' $\mathcal{O}_{\kappa}$, with
\[ c= c_{\kappa}=\frac{\kappa/k-1/2}{\kappa/k+1/2}. \]
\end{theorem}

To keep the article light and readable, we present our arguments in the context of Theorem~\ref{mainthm} first and then indicate the adjustments needed for the proof of Theorem~\ref{thinorbitthm} in Section~\ref{ThinOrbitsSection}.

\section{Preliminaries}

\subsection*{Notation} Throughout the paper, $\epsilon>0$ denotes a parameter which may be chosen to be as small as we like, but need not have the same value from one occurrence to another. The letter $p$ denotes an odd prime and $q=p^k$. We use $ \boldsymbol{\mu}_{p-1}$ to denote the set of $(p-1)^{\text{th}}$ roots of unity in the $p$-adic integers $\mathbb{Z}_p$. All implicit constants may depend on $\epsilon$ , $p$ and the parameter $\theta$ introduced in (\ref{moldef}), but not on $k$.

\subsection{Approximate functional equations}

We have the following standard approximate functional equations.

\begin{lemma}\label{afe}
 For a primitive Dirichlet character $\chi$ modulo $q$, let $\iota$ be defined by $\chi(-1)=(-1)^\iota$, and let
\begin{align}
\label{udef} &U(x)=\frac{1}{2\pi i} \int_{(2)} \frac{\Gamma(\frac{s+\iota}{2}+\frac{1}{4})}{\Gamma(\frac{\iota}{2}+\frac{1}{4})} (\pi^{\frac{1}{2}} x)^{-s} \frac{ds}{s},\\
\nonumber &V(x)=\frac{1}{2\pi i} \int_{(2)} \frac{\Gamma(\frac{s+\iota}{2}+\frac{1}{4})^2}{\Gamma(\frac{\iota}{2}+\frac{1}{4})^2} (\pi x)^{-s} \frac{ds}{s}.
\end{align}
We have that
\begin{align}
\label{u-v-estimates} U(x)\ll_{c} x^{-c}, \ \ \ \ V(x)\ll_{c} x^{-c}
\end{align}
for any $x,c>0$. For any $\lambda>0$, we have that
\begin{align}
&\label{afe1} L(\thalf, \chi) = \sum_{n\ge 1} \frac{\chi(n)}{n^\half} U\Big(\frac{n}{q^{1+\lambda}}\Big)+O(q^{-100}),\\
&\label{afe2} |L(\thalf,\chi)|^2 = 2\sum_{n_1,n_2\ge 1} \frac{\chi(n_1)\overline{\chi}(n_2)}{(n_1n_2)^\half} V\Big(\frac{n_1n_2}{q}\Big),
\end{align}
where the implied constant in \eqref{afe1} depends on $\lambda$.
\end{lemma}
\proof
The first equation (\ref{afe1}) is established by the functional equation of $L(s,\chi)$, which may be found in \cite[Theorem 4.15]{iwakow}, together with \cite[Theorem 5.3]{iwakow}, in which we take $G(u)=1$ and $X=q^{\half +\lambda}$. With this choice of $X$, the second sum in \cite[(5.12)]{iwakow} may be bounded by $q^{-100}$. The second equation is established by applying \cite[Theorem 5.3]{iwakow} to the product $L(s,\chi)L(s,\overline{\chi})$, rather than to each factor individually, with $G(u)=1$ and $X=1$. The estimates (\ref{u-v-estimates}) may be found in \cite[Proposition 5.4]{iwakow}.
\endproof

\noindent The sums in (\ref{afe1}) and (\ref{afe2}) are essentially restricted to $n<q^{1+\lambda+\epsilon}$ and $n_1n_2<q^{1+\epsilon}$, by (\ref{u-v-estimates}).

\subsection{Character averages} In this section, we record the orthogonality relations provided by averaging over the family of Dirichlet characters in a Galois orbit; see Lemma~\ref{chiavg} below. We start with a familiar auxiliary result.

\begin{lemma} \label{refine}
Let $m$ and $k\ge 1$ be integers. If $p^k\mid (m^p - 1)$, then $p^{k-1}\mid (m-1)$.
\end{lemma}

\proof

The claim is trivially true for $m=1$, so assume $m>1$.
Since $m\equiv m^p\equiv 1\bmod p$, we may write $m=1+p^rt$ for some $r\geqslant 1$ and some integer $t$ with $p\nmid t$. Then we have
\[ m^p-1=(1+p^rt)^p-1\equiv p^{r+1}t \  \bmod p^{r+2} , \]
so that $p^{r+1}$ is the highest power of $p$ that divides $m^p-1$. In particular, since $p^k\mid (m^p-1)$, we get $r+1\geqslant k$, and so $r\geqslant k-1$ and $p^{k-1}\mid(m-1)$.

\endproof

\begin{lemma}\label{chiavg} Let $q=p^k$ for an odd prime $p$ and let $\mathcal{O}$ be any Galois orbit of primitive Dirichlet characters mod $q$. For any integer $n$, we have that
\begin{align}
\label{chi-avg-zero} \sum_{\chi\in \mathcal{O}}  \chi(n) =0
\end{align}
unless 
\begin{align}
\label{cong-strong} n^{p-1}\equiv 1 \bmod p^{k-1}.
\end{align}
\end{lemma}
\proof
Suppose that $\mathcal{O}$ has cardinality given by (\ref{size}), for some $d\mid (p-1)$. We first show that (\ref{chi-avg-zero}) holds unless
\begin{align}
\label{cong-weak} n^{p(p-1)}\equiv 1 \bmod q.
\end{align}
It was shown in \cite[pg 17]{chi} that
\begin{align}
\label{chi-avg-ident}\frac{1}{|\mathcal{O}|}  \sum_{\chi\in \mathcal{O}}  \chi(n) = \frac{\mu(\text{ord}(n^{\frac{p-1}{d}}))}{\phi(\text{ord}(n^{\frac{p-1}{d}}))},
\end{align}
where ord$(n^{\frac{p-1}{d}})$ denotes the multiplicative order of $n^{\frac{p-1}{d}}$ in the group $(\mathbb{Z}/q\mathbb{Z})^\times$. Since $|(\mathbb{Z}/q\mathbb{Z})^\times|=p^{k-1}(p-1)$, if (\ref{cong-weak}) is not satisfied, then $p^2$ divides ord$(n^{\frac{p-1}{d}})$, and so (\ref{chi-avg-ident}) is zero.

Now if (\ref{cong-weak}) holds, then Lemma \ref{refine} implies, by taking $m=n^{p-1}$, that (\ref{cong-strong}) holds.
\endproof

Note that the condition \eqref{cong-strong}, which is all we will need from our orthogonality relations, requires less information than what is provided by the full sum over all $\chi\in\mathcal{O}$ or the explicit evaluation \eqref{chi-avg-ident}. This will be transparent in Lemma~\ref{chiavgthin}, which features a thinner average and for which we provide an independent proof.

\subsection{Mollifiers}
\label{MollifierSubsection}

The starting point of the mollifier method is the observation that
\begin{align}
\label{ratio} \frac{1}{|\mathcal{O}|}   \sum_{\substack{\chi\in \mathcal{O} \\ L(\half, \chi)\neq 0}} 1 \ \ge  \frac{1}{|\mathcal{O}|}   \sum_{\substack{\chi\in \mathcal{O} \\ L(\half, \chi) M(\half,\chi) \neq 0}} 1 \  \ge \frac{\Big|  \frac{1}{|\mathcal{O}|}   \sum_{ \chi\in \mathcal{O} } L(\half, \chi)M(\chi) \Big|^2}{  \frac{1}{|\mathcal{O}|}   \sum_{ \chi\in \mathcal{O} } | L(\half, \chi)|^2|M(\chi)|^2}.
\end{align} 
The second inequality above is the Cauchy-Schwarz inequality. The next step of the method is to evaluate the mollified moments. We prove the following. 
\begin{proposition}\label{1stmol} Let $q=p^k$ for an odd prime $p$ and let $\mathcal{O}$ be any Galois orbit of primitive Dirichlet characters mod $q$. For $0\le \theta<1$ in \eqref{moldef}, we have that
\begin{align*}
 \frac{1}{|\mathcal{O}|}  \sum_{\chi\in \mathcal{O}} L(\thalf, \chi) M(\chi) = 1 + O\big(q^{-\frac14+\epsilon} + q^{-\half +\frac{\theta}{2}+\epsilon}  \big).
 \end{align*}
\end{proposition}
\begin{proposition}\label{2ndmol} Let $q=p^k$ for an odd prime $p$ and let $\mathcal{O}$ be any Galois orbit of primitive Dirichlet characters mod $q$. Suppose that $\chi(-1)=(-1)^\iota$ for any $\chi\in \mathcal{O}$. For $0\le \theta<\half$ in \eqref{moldef}, we have that
\begin{multline*}
 \frac{1}{|\mathcal{O}|}  \sum_{\chi\in \mathcal{O}} |L(\thalf, \chi)|^2 |M(\chi)|^2 = \frac{p-1}{p}  \sum_{\substack{m_1,m_2\leq q^\theta\\(m_1m_2,q)=1}}  \frac{a_{m_1}\overline{a}_{m_2}}{[m_1,m_2]}\Big( \log\Big(\frac{q(m_1,m_2)^2}{\pi m_1m_2}\Big)+C\Big) \\+O\big(q^{-\frac14+\epsilon} +q^{-\half+\theta+\epsilon}\big),
 \end{multline*}
 where $(m_1,m_2)$ denotes the greatest common divisor of $m_1$ and $m_2$, $[m_1,m_2]=\frac{m_1m_2}{(m_1,m_2)}$ denotes the least common multiple of $m_1$ and $m_2$, 
 and $C=\frac{\Gamma'(\frac{1+2\iota}{4})}{\Gamma(\frac{1+2\iota}{4})} +2\gamma+2\frac{\log p}{p-1}$.
\end{proposition}
\noindent The final step is to insert the main terms of Propositions \ref{1stmol} and \ref{2ndmol} into the ratio on the right hand side of (\ref{ratio}), and then choose the coefficients $a_m$ so that the ratio is maximized. We note that the main terms of our mollified moments are identical to those in Iwaniec and Sarnak's problem (see \cite[(5.5)]{iwasar}) because in both problems the main terms arise from the diagonal contributions. That is, if $\chi(-1)=(-1)^{\iota}$ for all $\chi\in \mathcal{O}$, we have for any $0\le \theta < \half$ and some $\delta>0$ depending on $p$ and $\theta$ that
\begin{align*}
&\frac{1}{|\mathcal{O}|}  \sum_{\chi\in \mathcal{O}} L(\thalf, \chi) M(\chi) - \frac{2}{\varphi^{\star}(q)}  \summ_{\substack{ \chi \bmod q \\ \chi(-1)=(-1)^\iota}} L(\thalf, \chi) M(\chi) \ll q^{-\delta},\\
&\frac{1}{|\mathcal{O}|}  \sum_{\chi\in \mathcal{O}} |L(\thalf, \chi)|^2 |M(\chi)|^2 - \frac{2}{\varphi^{\star}(q)}  \summ_{\substack{ \chi \bmod q \\ \chi(-1)=(-1)^\iota}} |L(\thalf, \chi)|^2 |M(\chi)|^2 \ll q^{-\delta},
\end{align*}
 where $\summ$ means that summation is restricted to the primitive characters and $\varphi^{\star}(q)= \summ_{\chi \bmod q} 1$. The optimal choice for $a_m$ will therefore be the same as in Iwaniec and Sarnak's problem (and we will not repeat the proof of this here), which yields that the ratio on the right hand side of (\ref{ratio}) is $\frac{\theta}{1+\theta}+o(1)$ as $k\to\infty$. Taking $\theta\to\frac12$, Theorem \ref{mainthm} thus follows from Propositions \ref{1stmol} and \ref{2ndmol}.

\subsection{Roth's theorem}
\label{RothSection}

The following lemma is a special case of a result of Rohrlich \cite[Proposition 1, pg 401]{roh2}. It is a consequence of the $p$-adic version of Roth's theorem, due to Ridout \cite{rid}. 

\begin{lemma}\label{rohrlich-roth}
Let $\beta\in \mathbb{Z}_p$ and $0<\delta<\half$. Suppose that $\beta$ is algebraic over $\mathbb{Q}$ of degree at least 2. Then, for sufficiently large $k\ge k_0(\beta,\delta)$, there are no nonzero integers $a$ and $b$ which satisfy
\begin{align*}
|a -  b \beta |_p \le p^{-k+1}
\end{align*}
and
\begin{align*}
|a|,|b|< (p^{k})^{\half-\delta}.
\end{align*}
\end{lemma}
\proof
In \cite[Proposition 1, pg 401]{roh2}, take $\alpha_p=1$ and $\beta_p=\beta$. Note that condition ${\rm (ii)}$ directly above the proposition is satisfied because $\beta \notin \mathbb{Q}.$ Note that $k_0$ is not computable.
\endproof 
\noindent When $(b,p)=1$, this result says that an approximation of the algebraic number $\beta$ by $\frac{a}{b}$ to an error within $(\max\{|a|,|b|\})^{-2-\delta}$ is too good to exist. Note that there are infinitely many approximations of any $\beta\in\mathbb{Z}_p$ up to within $\ll(\max\{|a|,|b|\})^{-2}$ by a $p$-adic incarnation of Dirichlet's Approximation Theorem (see also~\cite{BlMi} for an analogue of the Farey dissection in this context), and it is a hallmark of Roth's theorem that the exponent in Lemma~\ref{rohrlich-roth} is essentially the best possible. For an interesting investigation into existence of algebraic numbers exhibiting approximability by rationals within the transition range between the exponents $-2$ and $-2-\delta$, 
see \cite{lag}.

We stress that the $p$-adic analogues of previous partial results toward Roth's Theorem, such as the theorems of Liouville, Thue, Siegel, and Dyson, do not (except for a few small primes $p$) suffice to obtain an asymptotic in the situation of Proposition~\ref{2ndmol} or Theorem~\ref{2ndthm} with our methods. This is so because, for example when estimating \eqref{subdyadic2}, in order to obtain \eqref{finaleq} with a $o(1)$-upper bound on the right-hand side for large $q$, we need to be able to take $Q$ (essentially the allowable length of intervals in Lemma~\ref{ridoutapplication}, below) to be at least $(p^k)^{\frac14+\delta}$, whereas all results prior to the actual $p$-adic Roth's theorem furnish only $(p^k)^{o_p(1)}$.

We will also use the following convenient implication of Lemma~\ref{rohrlich-roth}.

\begin{lemma} \label{ridoutapplication}
Let $0<\delta<\half$, let $k\ge k_1(\delta)$ be sufficiently large, and let $\mathcal{A}_k$ and $\mathcal{B}_k$ be intervals in the rational integers of length at most $(p^k)^{\half-\delta}$. Then there are at most $p-3$ pairs $(a,b)\in\mathcal{A}_k\times\mathcal{B}_k$ such that 
\begin{align*}
(ab,p)=1,
\end{align*}
\begin{align}
\label{notequiv} a\not\equiv\pm b \bmod p^{k-1},
\end{align}
and 
\begin{align}
\label{congruence} a^{p-1}- b^{p-1}\equiv 0 \bmod p^{k-1}.
\end{align}
\end{lemma}
\proof
Suppose that $a$ and $b$ satisfy the conditions of the lemma. The congruence (\ref{congruence}) implies that the $p$-adic integer $a^{p-1}- b^{p-1}$ has norm
\begin{align}
\label{factorize} \big|a^{p-1}- b^{p-1}\big|_p = \Big| \prod_{\zeta\in \boldsymbol{\mu}_{p-1}} (a-b\zeta) \Big|_p \leq p^{-k+1}.
\end{align}
Recall that the roots of unity in $\boldsymbol{\mu}_{p-1}$ are distinct modulo $p$. By this fact and the assumption $(ab,p)=1$, we have that $a-b\zeta$ and $a-b\zeta'$ are distinct modulo $p$ for $\zeta\neq \zeta'$. Thus, by (\ref{factorize}) and (\ref{notequiv}), for some $\zeta\in \boldsymbol{\mu}_{p-1}\backslash \{\pm 1\}$ we have
\begin{align*}
|a-b\zeta| _p \leq p^{-k+1}.
\end{align*}

Now suppose for a contradiction that at least $p-2$ pairs $(a,b)$ satisfy the conditions of the lemma. Then by the argument above and Dirichlet's Box Principle, we have for some $\zeta\in \boldsymbol{\mu}_{p-1}\backslash\{\pm 1\}$, at least two distinct pairs $(a_1,b_1),(a_2,b_2)\in \mathcal{A}_k\times \mathcal{B}_k$ satisfying
\begin{align*}
 |a_j-b_j\zeta |_p \leq p^{-k+1}  \qquad (j=1,2).
\end{align*}
It follows from the strong triangle inequality that
\begin{align*}
|(a_2-a_1)-(b_2-b_1)\zeta) |_p \le \max\{ |a_1-b_1\zeta |_p, |a_2-b_2\zeta |_p \} \leq p^{-k+1}.
\end{align*}
We also have, by the assumption on the lengths of $\mathcal{A}_k$ and $\mathcal{B}_k$, that
\begin{align*}
|a_2-a_1|, |b_2-b_1|\le (p^k)^{\frac{1}{2}-\delta}.
\end{align*}
Thus by Lemma \ref{rohrlich-roth}, in which we take $\beta=\zeta$, we deduce that for $k>k_1$ sufficiently large, where $k_1$ is not computable, we must have $a_2-a_1=b_2-b_1=0$. This is the desired contradiction.
\endproof

\section{Proofs of the main results}

\subsection{Proof of Proposition \ref{1stmol}}

In Lemma \ref{afe}, let $\lambda>0$ be such that $\lambda+\theta<1$. We have that
\begin{align}
\label{1stmol-expansion} \frac{1}{|\mathcal{O}|}  \sum_{\chi\in \mathcal{O}} L(\thalf, \chi) M(\chi) = \sum_{m\leq q^\theta,} \sum_{n\ge 1} \frac{a_m}{(nm)^\half} U\Big(\frac{n}{q^{1+\lambda}}\Big) \frac{1}{|\mathcal{O}|}  \sum_{\chi\in \mathcal{O}} \chi(nm) + O(q^{-99}).
\end{align}
Note that the exchange of summation above is valid because although the function $U(x)$ depends on $\chi(-1)$, this value is the same for every character in $\mathcal{O}$. Now, by Lemma \ref{chiavg}, we have that 
\begin{align*}
 \frac{1}{|\mathcal{O}|} \sum_{\chi\in \mathcal{O}} \chi(nm) = 0 
\end{align*}
unless $(nm)^{p-1}\equiv 1\bmod{p^{k-1}}$, in which case we have by Hensel's Lemma that $nm$ is congruent to one of $\zeta_1,\ldots, \zeta_{p-1}\in\boldsymbol{\mu}_{p-1}$ modulo $p^{k-1}\mathbb{Z}$. Separating the term $nm=1$ for the main term in Proposition~\ref{1stmol}, and using also the divisor bound, we have that (\ref{1stmol-expansion}) equals
\begin{equation}
\label{summed1st}
U\Big(\frac{1}{q^{1+\lambda}}\Big)
+O\Big(q^{\epsilon}\sum_{1\leq j\leq p-1}\sum_{\substack{1<s\leq q\\ s\equiv\zeta_j\bmod p^{k-1}}}s^{-\half}\Big)
+O\Big(q^{\epsilon} \sum_{1\le j\le p-1} \sum_{\substack{q<s<q^{1+\lambda+\theta+\epsilon}\\ s\equiv \zeta_j \bmod p^{k-1}}} s^{-\half} \Big).
\end{equation}
By shifting the line of integration in (\ref{udef}) to $\Re(s)=-\half+\epsilon$, we see that the main term above equals $1+O(q^{-\half+\epsilon})$.

The summands in the first error term are particularly sensitive to the size of $s$. From the condition that $s^{p-1}\equiv 1\bmod{p^{k-1}}$ and $s>1$, we immediately have that $s>p^{\frac{k-1}{p-1}}$, and the first error term is seen to be
\[ O\Big(q^{-\frac1{2(p-1)}+\epsilon}\Big), \]
without recourse to $p$-adic Roth's theorem. However, we can improve this estimate by appealing to Lemma~\ref{rohrlich-roth}. The terms with $s\equiv\pm 1\bmod{p^{k-1}}$, $s>1$, contribute $O(p^{-(k-1)/2})$. As for the terms corresponding to $s\equiv\zeta_j\not\equiv\pm 1\bmod{p^{k-1}}$, Lemma~\ref{rohrlich-roth} with $\beta=\zeta_j$ guarantees that, for sufficiently large $k$, there are no values of $|s|<(p^k)^{\frac12-\delta}$ with $|s-\zeta_j\cdot 1|_p\leq p^{-k+1}$, so that all these terms must in fact satisfy $s\geq (p^k)^{\frac12-\delta}$, and in total the first error term is
\[ O\Big(q^{-\frac14+\frac12\delta+\epsilon}\Big). \]

Finally, writing $s=\zeta_j+p^{k-1}r$, we have that the second error term above is bounded by
\begin{align*}
q^{\epsilon} \sum_{1\le r < q^{\lambda + \theta+\epsilon} } (qr)^{-\half} \ll q^{\frac{\lambda+\theta-1}{2}+\epsilon}.
\end{align*}
Taking the positive $\lambda$ and $\delta$ to be as small as we like and adjusting the implied constants completes the proof. 
\qed

\subsection{Proof of Proposition \ref{2ndmol}}

By Lemma \ref{afe}, we have that
\begin{equation}
\label{funceq2ndmol}
\begin{aligned}
 \frac{1}{|\mathcal{O}|}  \sum_{\chi\in \mathcal{O}} &|L(\thalf, \chi)|^2 |M(\chi)|^2 \\&{}= 2 \sum_{m_1,m_2 \leq q^\theta,} \sum_{n_1,n_2\ge 1} \frac{a_{m_1}\overline{a}_{m_2}}{(n_1n_2m_1m_2)^\half}V\Big(\frac{n_1n_2}{q}\Big)  \frac{1}{|\mathcal{O}|}  \sum_{\chi\in \mathcal{O}} \chi(n_1m_1)\overline{\chi}(n_2m_2).
\end{aligned}
\end{equation}
We write this as a sum of diagonal terms and off-diagonal terms,
\begin{align*}
\sum_{n_1m_1=n_2m_2}+\sum_{n_1m_1 \neq n_2m_2}.
\end{align*}
\subsubsection*{\bf The diagonal} We first consider the diagonal terms. The equality $n_1m_1=n_2m_2$ is the same as requiring $n_1=r m_2/(m_1,m_2)$ and $n_2=r m_1/(m_1,m_2)$ for some $r\in\mathbb{N}$. Thus the diagonal contribution is
\begin{align*}
\sum_{n_1m_1=n_2m_2}= 2\sum_{\substack{m_1,m_2 \leq q^\theta\\ (m_1m_2,p)=1}}   \frac{a_{m_1}\overline{a}_{m_2} }{[m_1,m_2]} \sum_{\substack{r\ge 1\\ (r,p)=1}} \frac{1}{r}V\Big( \frac{r^2 m_1m_2}{q(m_1,m_2)^2} \Big).
\end{align*}
The innermost sum above can be evaluated by the calculation in \cite[Lemma 4.1]{iwasar} and equals 
\begin{align*}
\frac{p-1}{2p} \Big( \log\Big(\frac{q(m_1,m_2)^2}{\pi m_1m_2}\Big)+C\Big) +O\Big(\Big(\frac{q(m_1,m_2)^2}{m_1m_2}\Big)^{-\frac{1}{2}+\epsilon}\Big),
\end{align*}
where $C$ is as in the statement of Proposition \ref{2ndmol}. The main term above gives the main term of Proposition \ref{2ndmol}, and the total error is less than
\begin{align*}
q^\epsilon \sum_{m_1,m_2 \leq q^\theta} \frac{1}{[m_1,m_2]}\Big( \frac{q(m_1,m_2)^2}{m_1m_2}\Big)^{-\frac{1}{2}}\ll q^{-\half+\epsilon} \sum_{m_1,m_2 \leq q^\theta} \frac{1}{(m_1m_2)^{\half}}\ll q^{-\half+\theta+\epsilon}.
\end{align*}

\subsubsection*{\bf The off-diagonal} We now turn to the off-diagonal contribution, which we must bound by a negative power of $q$. By Lemma \ref{chiavg} and \eqref{u-v-estimates}, we have that 
\begin{align*}
\sum_{n_1m_1 \neq n_2m_2}\ll  \mathop{q^\epsilon \sum_{\substack{m_1,m_2 \leq q^\theta,}}   \sum_{\substack{n_1n_2\leq q^{1+\epsilon}}} }_{\substack{(n_1n_2m_1m_2,p)=1\\ (n_1m_1)^{p-1}\equiv (n_2m_2)^{p-1} \bmod p^{k-1} \\ n_1m_1\neq n_2m_2}} \frac{1}{(n_1n_2m_1m_2)^\half}.
\end{align*}
Writing $a=n_1m_1$ and $b=n_2m_2$ and splitting the resulting sum above into dyadic intervals $A\le a < 2A$ and $B\le b < 2 B$, it suffices to show that the sum
\begin{align}
\label{dyadic0} \frac{q^\epsilon}{(AB)^\half}\sum_{\substack{A\le a < 2A\\   B\le b < 2B\\ (ab,p)=1\\ a^{p-1}\equiv b^{p-1}\bmod p^{k-1}\\ a\neq b}} 1
\end{align}
is less than a negative power of $q$ for 
\begin{align}
\label{restriction} 1\le AB \le q^{1+2\theta+\epsilon}.
\end{align}
We split (\ref{dyadic0}) further as
\begin{align}
\label{dyadic-split} \frac{q^\epsilon}{(AB)^\half}\sum_{\substack{A\le a < 2A\\   B\le b < 2B\\ (ab,p)=1\\ a \equiv\pm b \bmod p^{k-1}\\ a\neq b}} 1 \indent  + \indent  \frac{q^\epsilon}{(AB)^\half}\sum_{\substack{A\le a < 2A\\   B\le b < 2B\\ (ab,p)=1\\ a^{p-1}\equiv b^{p-1}\bmod p^{k-1}\\ a \not\equiv\pm b \bmod p^{k-1} }} 1.
\end{align}

In the first sum above, we must have $2A>p^{k-1}$ or $2B>p^{k-1}$.  In the first case, for each of the $B$ choices for $b$, there are $O(\frac{A}{q})$ choices for $a$ by the Chinese Remainder Theorem. In the same way in the second case there are $O(\frac{AB}{q})$ choices for $a$ and $b$. Thus the first sum in (\ref{dyadic-split}) is bounded by
\begin{align*}
\frac{q^\epsilon}{(AB)^\half} \frac{AB}{q} \ll q^{-\half+\theta+\epsilon}.
\end{align*}

Now we consider the second sum in (\ref{dyadic-split}). We must show that
\begin{align}
\label{dyadic} \frac{q^\epsilon}{(AB)^\half}\sum_{\substack{A\le a < 2A\\   B\le b < 2B\\ (ab,p)=1\\ a^{p-1}\equiv b^{p-1}\bmod p^{k-1}\\ a \not\equiv\pm b \bmod p^{k-1} }} 1
\end{align}
is bounded by a negative power of $q$. 
We may assume that 
\begin{align}
\label{assume} q^{-\half} <\frac{A}{B} < q^\half,
\end{align}
since otherwise the proof is complete. To see this, suppose without loss of generality that $A\le B$. Then for each of the $A$ choices of $a$ in (\ref{dyadic}), there are $O(1+\frac{B}{q})$ possibilities of $b$ which satisfy the congruence $a^{p-1}\equiv b^{p-1}\bmod{p^{k-1}}$.
Thus if (\ref{assume}) is not satisfied, then (\ref{dyadic}) is bounded by 
\begin{equation}
\label{mayassume}
 \frac{q^\epsilon}{(AB)^\half} \Big(A+\frac{AB}{q}\Big)  \ll q^{-\frac{1}{4}+\epsilon}+ q^{-\half+\theta+\epsilon}.
\end{equation}

Now, 
assuming (\ref{assume}), we analyze (\ref{dyadic}) according to the sizes of $A$ and $B$ as follows. Let $0<\delta<\frac{1}{p}$ and $Q=q^{\frac{1}{2}-2\delta}$.

\subsubsection*{Case $1$} Suppose that $2A<q^{\half-\delta}$ and $2B<q^{\half-\delta}$. Every term in \eqref{dyadic} would have to satisfy $a\equiv b\zeta_j\bmod{p^{k-1}}$ for some $\zeta_j\in\boldsymbol{\mu}_{p-1}\setminus\{\pm 1\}$, and hence $|a-b\zeta_j|_p\leq p^{-k+1}$ with $|a|,|b|< (p^k)^{\frac12-\delta}$. For sufficiently large $k$, however, by Lemma \ref{rohrlich-roth} there are no such nonzero integers $a$ and $b$, and the corresponding sum \eqref{dyadic} is actually empty.

\subsubsection*{Case $2$} Suppose that $2A\ge q^{\half-\delta}$ and $2B< q^{\half-\delta}$. (The case $2A<q^{\half-\delta}$ and $2B\ge q^{\half-\delta}$ is treated similarly). Dividing the dyadic interval $A\le a < 2A$ into smaller pieces of length $Q$, we may bound (\ref{dyadic}) by 
\begin{align}
\label{subdyadic1} \frac{q^\epsilon}{(AB)^\half} \sum_{\substack{ 1\le u \le \frac{A}{Q},}} \sum_{\substack{A+(u-1)Q\le a < A+uQ\\   B\le b < 2B \\ (ab,p)=1\\ a^{p-1}\equiv b^{p-1}\bmod p^{k-1} \\ a \not\equiv\pm b \bmod p^{k-1}}} 1.
\end{align}
By Lemma \ref{ridoutapplication}, the innermost sum of (\ref{subdyadic1}) is, for sufficiently large $k$, bounded by a constant (depending on $p$). Using this fact and (\ref{assume}), we see that (\ref{subdyadic1}) is bounded by
\begin{align}
\label{almostfinaleq}
  \frac{q^\epsilon}{(AB)^\half} \frac{A}{Q} \ll q^{-\frac{1}{4}+2\delta+\epsilon}.
\end{align}
This falls into the error term of Proposition \ref{2ndmol} as we may take $\delta$ to be as small as we like.

\subsubsection*{Case $3$} Suppose that $2A\ge q^{\half-\delta}$ and $2B\ge q^{\half-\delta}$. Dividing the dyadic intervals $A\le a < 2A$ and $B\le b < 2B$ into smaller pieces of length $Q$, we may rewrite (\ref{dyadic}) as 
\begin{align}
\label{subdyadic2} \frac{q^\epsilon}{(AB)^\half} \sum_{\substack{ 1\le u \le \frac{A}{Q},\\ 1\le v\le \frac{B}{Q},}} \sum_{\substack{A+(u-1)Q\le a < A+uQ\\   B+(v-1)Q\le b < B+vQ\\ (ab,p)=1\\ a^{p-1}\equiv b^{p-1}\bmod p^{k-1} \\ a \not\equiv\pm b \bmod p^{k-1} }} 1.
\end{align}
By Lemma \ref{ridoutapplication}, the innermost sum of (\ref{subdyadic1}) is, for sufficiently large $k$, bounded by a constant (depending on $p$). Using this fact and (\ref{restriction}), we see that (\ref{subdyadic2}) is bounded by
\begin{align}
\label{finaleq}
\frac{q^\epsilon}{(AB)^\half} \frac{AB}{Q^2} \ll q^{-\half+\theta+4\delta+\epsilon}.
\end{align}
This falls into the error term of Proposition \ref{2ndmol} as we may take $\delta$ to be as small as we like.

The proof of Proposition \ref{2ndmol} is now complete.
\qed

\subsection{The case of thin orbits}
\label{ThinOrbitsSection}

In this section, we prove Theorem~\ref{thinorbitthm}. As is customary in analytic number theory (and as is already the case with full Galois orbits of primitive characters), the principal change introduced by the shrinking family of characters in the orbit $\mathcal{O}_{\kappa}$ is that more terms survive averaging over the family. We quantify this effect with the following modification of Lemma~\ref{chiavg} on orthogonality relations.

\begin{lemma}\label{chiavgthin} Let $q=p^k$ for an odd prime $p$, let $0<\kappa\leq k-1$, and let $\mathcal{O}_{\kappa}$ be a thin Galois orbit of primitive Dirichlet characters mod $q$. For any integer $n$, we have that
\begin{align}
\label{chi-avg-zero-thin} \sum_{\chi\in \mathcal{O}_{\kappa}}  \chi(n) =0
\end{align}
unless 
\begin{align}
\label{cong-strong-thin} n^{p-1}\equiv 1 \bmod p^{\tilde{\kappa}+1},
\end{align}
where $\tilde{\kappa}=\min(\kappa,k-2)$.
\end{lemma}
\proof

Fix a character $\chi_0\in \mathcal{O}_{\kappa}$. 
Fix a generator $g$ of the cyclic group $(\mathbb{Z}/p^k\mathbb{Z})^\times$, and write $\chi_0(g)=\xi^{\gamma}$ for some $\gamma$. In particular, $(\gamma,\phi(q))=(p-1)/d$, where $d\mid (p-1)$ is the order of $\chi_0$ and of all characters in $\mathcal{O}_{\kappa}$, so that the corresponding full orbit $\mathcal{O}\supseteq\mathcal{O}_{\kappa}$ has cardinality given by \eqref{size}.

If $n$ is an integer divisible by $p$, the lemma is trivially true. We therefore assume $n$ is relatively prime to $p$ and let $0\leq r< p^{k-1}(p-1)$ be such that $n=g^r$ in $(\mathbb{Z}/p^k\mathbb{Z})^\times$. 
From \eqref{SigmaXiCondition}, it is immediate that
 \[ \sum_{\chi\in\mathcal{O}_{\kappa}}\chi(n)=
 \sum_{\substack{a\bmod p^{k-1}(p-1)\\ a\equiv 1\bmod{p^{k-1-\kappa}(p-1)}}}\xi^{\gamma ra}=
\begin{cases}
	\chi_0(n) 
	\sum\limits_{0\leq j<p^\kappa} e\big(\gamma rj/p^{\kappa}\big), &0\leqslant\kappa<k-1,\\
	\chi_0(n)^p\sum\limits_{0\leq j<p^{k-1},\,p\nmid j}e\big(\gamma rj/p^{k-1}\big), &\kappa=k-1.
\end{cases}	
	\]
In either case, the resulting sum vanishes unless $p^{\tilde{\kappa}}\mid r$. The condition $p^{\tilde{\kappa}} \mid r$ is equivalent to $n^{p^{k-1-\tilde{\kappa}}(p-1)}\equiv 1 \bmod p^k$, which by Lemma \ref{refine} implies that 
$n^{p-1}\equiv 1 \bmod p^{\tilde{\kappa}+1}$. This proves the lemma.
\endproof

We will only use \eqref{cong-strong-thin} as a condition modulo $p^{\kappa}$. In particular, note that the localization \eqref{cong-strong} is achieved already by averaging over $\chi\in\mathcal{O}_{k-1}$. The resulting Ramanujan sum in the case $\kappa=k-1$ should be compared with the explicit evaluation \eqref{chi-avg-ident}.

The analogs of Propositions~\ref{1stmol} and \ref{2ndmol} on mollified moments are as follows:

\begin{proposition}\label{thinmol} Let $q=p^k$ for an odd prime $p$, let $k/2<\kappa\leq k-1$, and let $\mathcal{O}_{\kappa}$ be a thin Galois orbit of primitive Dirichlet characters mod $q$. For $0\le \theta<2(\frac{\kappa}k-\half)$ in \eqref{moldef}, we have
\begin{align*}
 \frac{1}{|\mathcal{O}_{\kappa}|}  \sum_{\chi\in \mathcal{O}_{\kappa}} L(\thalf, \chi) M(\chi) = 1 + O\big(q^{-\frac{\kappa}{4k}+\epsilon} + q^{\half+\frac{\theta}{2}-\frac{\kappa}k+\epsilon}  \big),
 \end{align*}
while, for $0\le \theta<\frac{\kappa}k-\half$ in \eqref{moldef}, we have
\begin{multline*}
 \frac{1}{|\mathcal{O}_{\kappa}|}  \sum_{\chi\in \mathcal{O}_{\kappa}} |L(\thalf, \chi)|^2 |M(\chi)|^2 = \frac{p-1}{p} \sum_{\substack{m_1,m_2\leq q^\theta\\(m_1m_2,q)=1}}  \frac{a_{m_1}\overline{a}_{m_2}}{[m_1,m_2]}\Big( \log\Big(\frac{q(m_1,m_2)^2}{\pi m_1m_2}\Big)+C\Big)\\+O\big(q^{-\frac{\kappa}{4k}+\epsilon} +q^{\half+\theta-\frac{\kappa}k+\epsilon}\big),
 \end{multline*}
with notations as in Proposition~\ref{2ndmol}.
\end{proposition}

\proof
The proof follows the proofs of Propositions~\ref{1stmol} and \ref{2ndmol}, with Lemma~\ref{chiavgthin} as the orthogonality relation in place of Lemma~\ref{chiavg}.

For the first mollified moment, we start by inserting the approximate functional equation as in \eqref{1stmol-expansion}. By Lemma~\ref{chiavgthin}, the character average isolates the main term, which comes from $mn=1$ and is identical as before, and summands with $s=mn\equiv\zeta_j\bmod{p^{\kappa}}$, $s>1$, which we split into two terms, corresponding to the ranges $1<s\leq p^{\kappa}$ and $p^{\kappa}<s<q^{1+\lambda+\theta+\epsilon}$. By Lemma~\ref{rohrlich-roth}, for sufficiently large $k$, all summands in the first sum satisfy $s\geq (p^{\kappa})^{\half-\delta}$, and in total the first error term is $O\big((p^{\kappa})^{-\frac14+\frac12\delta}\big)$. Writing $s=\zeta_j+p^{\kappa}r$, the second error term is similarly bounded by
\[ q^{\epsilon}\sum_{1\leq r<q^{1+\lambda+\theta+\epsilon-\kappa/k}}(p^{\kappa}r)^{-\half}\ll q^{\frac{\lambda+\theta+1}2-\frac{\kappa}k}, \]
completing the proof of the first part of our proposition.

As for the second mollified moment, we start by expanding using the same functional equation~\eqref{funceq2ndmol}. The main term arises from the diagonal terms, when $n_1m_1=n_2m_2$, and is the same as above, while in the off-diagonal terms we are left by Lemma~\ref{chiavgthin} with the sums as in \eqref{dyadic0} now subject to $a^{p-1}\equiv b^{p-1}\bmod{p^{\kappa}}$, and which we further split as in \eqref{dyadic-split}. In the first sum, over $a\equiv\pm b\bmod{p^{\kappa}}$, we must have $2A>p^{\kappa}$ or $2B>p^{\kappa}$, and the total contribution from these terms is bounded by
\[ \frac{q^{\epsilon}}{(AB)^{\half}}\frac{AB}{p^{\kappa}}\ll q^{\half+\theta-\frac{\kappa}k+\epsilon}. \]
In the second sum, analogously as before, we may assume that $p^{-\frac{\kappa}2}<A/B<p^{\frac{\kappa}2}$, since otherwise the proof is complete as in \eqref{mayassume}. We now set $Q=(p^{\kappa})^{\frac12-2\delta}$ and consider three cases. If $2A<(p^{\kappa})^{\half-\delta}$ and $2B<(p^{\kappa})^{\half-\delta}$, then the sum is empty for sufficiently large $k$ by Lemma~\ref{rohrlich-roth}. If $2A\geq (p^{\kappa})^{\half-\delta}$ and $2B<(p^{\kappa})^{\half-\delta}$, then by splitting the $a$-sum into intervals of length $Q$ the sum is bounded as in \eqref{almostfinaleq} by
\[ \frac{q^{\epsilon}}{(AB)^{\half}}\frac AQ\ll q^{-\frac{\kappa}{4k}+2\delta+\epsilon}, \]
and, finally, if both $2A\geq (p^{\kappa})^{\half-\delta}$ and $2B\geq (p^{\kappa})^{\half-\delta}$, then the sum is bounded as in \eqref{finaleq} by
\[ \frac{q^{\epsilon}}{(AB)^{\half}}\frac{AB}{Q^2}\ll q^{\half+\theta-\frac{\kappa}k+4\delta\frac{\kappa}k+\epsilon}, \]
which completes the proof as we may take $\delta$ as small as we wish.
\endproof

Theorem~\ref{thinorbitthm} now follows from Proposition~\ref{thinmol} in the same way in which Theorem~\ref{mainthm} is deduced from Propositions~\ref{1stmol} and \ref{2ndmol}, by taking $\theta\to\frac{\kappa}k-\half$.

\subsection*{Acknowledgements} Part of this work was done during the excellent National Science Foundation-supported $29^{\rm th}$ Automorphic Forms Workshop held at the University of Michigan, Ann Arbor; we wish to express our thanks to the organizers of the conference. The third author is grateful to Jeffrey C. Lagarias for his constant encouragement. 

\bibliographystyle{amsplain}

\bibliography{galois-orbits}

\end{document}